\documentclass[11pt]{article}
\usepackage{amsmath,amssymb,amsthm,amscd}

\newtheorem{thm}{Theorem}[section]
\newtheorem{lem}[thm]{Lemma}

\theoremstyle{definition}
\newtheorem{df}[thm]{Definition}
\newtheorem{rem}[thm]{Remark}

\numberwithin{equation}{section}

\renewcommand{\phi}{\varphi}

\newcommand{\ep}{\varepsilon}

\newcommand{\Ad}{\operatorname{Ad}}
\newcommand{\Aut}{\operatorname{Aut}}
\newcommand{\Ext}{\operatorname{Ext}}
\newcommand{\id}{\operatorname{id}}
\newcommand{\Hom}{\operatorname{Hom}}
\newcommand{\Ima}{\operatorname{Im}}
\newcommand{\Ker}{\operatorname{Ker}}
\newcommand{\Lip}{\operatorname{Lip}}
\newcommand{\OrderExt}{\operatorname{OrderExt}}
\newcommand{\Sp}{\operatorname{Sp}}

\newcommand{\N}{\mathbb{N}}
\newcommand{\Z}{\mathbb{Z}}
\newcommand{\Q}{\mathbb{Q}}
\newcommand{\R}{\mathbb{R}}
\newcommand{\C}{\mathbb{C}}
\newcommand{\T}{\mathbb{T}}

\title{$\Z^N$-actions on UHF algebras of infinite type}
\author{Hiroki Matui \\
Graduate School of Science \\
Chiba University \\
Inage-ku, Chiba 263-8522, Japan}
\date{}

\begin{document}
\maketitle

\begin{abstract}
We prove that 
all strongly outer $\Z^N$-actions on a UHF algebra of infinite type are 
strongly cocycle conjugate to each other. 
We also prove that 
all strongly outer, asymptotically representable $\Z^N$-actions 
on a unital simple AH algebra with real rank zero, slow dimension growth 
and finitely many extremal tracial states are 
cocycle conjugate to each other. 
\end{abstract}

\section{Introduction}

Classification of group actions is 
one of the most fundamental subjects in the theory of operator algebras. 
For AFD factors, 
a complete classification is known for actions of countable amenable groups. 
However, classification of automorphisms or group actions 
on $C^*$-algebras is still a far less developed subject, 
partly because of $K$-theoretical difficulties. 
In this paper we prove the uniqueness of strongly outer $\Z^N$-actions 
on UHF algebras of infinite type (Theorem \ref{SCConUHF}) 
as well as the uniqueness of 
strongly outer, asymptotically representable $\Z^N$-actions 
on certain unital simple AH algebras (Theorem \ref{CCofasymp}). 

A. Kishimoto \cite{K95crelle} proved that 
any strongly outer $\Z$-actions on any UHF algebras have the Rohlin property 
and that they are strongly cocycle conjugate to each other. 
H. Nakamura \cite{N1} showed that 
any strongly outer $\Z^2$-actions on UHF algebras have the Rohlin property. 
T. Katsura and the author \cite{KM} gave a complete classification 
of strongly outer $\Z^2$-actions on UHF algebras 
by using the Rohlin property. 
In particular, it was shown that 
strongly outer $\Z^2$-actions on any UHF algebra of infinite type 
are unique up to cocycle conjugacy. 
The present paper generalizes these results. 
We show that strong outerness is equivalent to the Rohlin property 
for $\Z^N$-actions on UHF algebras of infinite type 
(Theorem \ref{Rohlintype1}) and 
that all strongly outer $\Z^N$-actions on any UHF algebra of infinite type 
are strongly cocycle conjugate to each other (Theorem \ref{SCConUHF}). 

We briefly review other classification results of $\Z^N$-actions known so far. 
For AT algebras, 
A. Kishimoto \cite{K98JOT,K98JFA} showed 
the Rohlin property for a certain class of automorphisms and 
obtained a cocycle conjugacy result. 
The author \cite{M10CMP} extended this result to 
unital simple AH algebras with real rank zero and slow dimension growth. 
A certain class of $\Z^2$-actions on unital simple AF algebras were 
also classified in \cite{M10CMP}. 
Y. Sato \cite{S} proved that 
strongly outer $\Z$-actions on the Jiang-Su algebra $\mathcal{Z}$ are 
unique up to strong cocycle conjugacy. 
Y. Sato and the author \cite{MS} obtained 
the uniqueness of strongly outer $\Z^2$-actions on $\mathcal{Z}$. 
For Kirchberg algebras, 
complete classification of aperiodic automorphisms was given 
by H. Nakamura \cite{N2}. 
M. Izumi and the author \cite{IM} classified 
a large class of $\Z^2$-actions 
and also showed the uniqueness of $\Z^N$-actions 
on $\mathcal{O}_2$, $\mathcal{O}_\infty$ and 
$\mathcal{O}_\infty\otimes B$ with $B$ being a UHF algebra of infinite type 
(see also \cite{M08}). 

This paper is organized as follows. 
In Section 2, we collect for reference 
basic notations, terminologies and definitions. 
In Section 3, it is shown that 
$\Z^N$-actions with the Rohlin property on a UHF algebra of infinite type 
are mutually cocycle conjugate (Theorem \ref{CCofRohlin1}). 
It is also shown that 
asymptotically representable $\Z^N$-actions with the Rohlin property 
on a unital simple AH algebra with real rank zero and slow dimension growth 
are mutually cocycle conjugate (Theorem \ref{CCofRohlin2}). 
In Section 4, we prove that 
strongly outer $\Z^N$-actions on UHF algebras of infinite type 
have the Rohlin property (Theorem \ref{Rohlintype1}). 
The same statement is also obtained 
for approximately representable, strongly outer $\Z^N$-actions 
on certain AH algebras (Theorem \ref{Rohlintype2}). 
This, together with the result of Section 3, implies 
the uniqueness of asymptotically representable, strongly outer $\Z^N$-actions 
on AH algebras (Theorem \ref{CCofasymp}). 
In Section 5, 
we show a kind of cohomology vanishing theorem (Lemma \ref{admissible}) and 
complete the proof of the uniqueness of strongly outer $\Z^N$-actions 
on UHF algebras of infinite type (Theorem \ref{SCConUHF}).

\section{Preliminaries}

The cardinality of a set $F$ is written by $\lvert F\rvert$. 
For a Lipschitz continuous function $f$, 
we denote its Lipschitz constant by $\Lip(f)$. 

Let $A$ be a $C^*$-algebra. 
For $a,b\in A$, we mean by $[a,b]$ the commutator $ab-ba$. 
The set of tracial states on $A$ is denoted by $T(A)$. 
When $A$ is unital, 
we mean by $U(A)$ the set of all unitaries of $A$. 
For $u\in U(A)$, 
the inner automorphism induced by $u$ is written by $\Ad u$. 
An automorphism $\alpha\in\Aut(A)$ is called outer, 
when it is not inner. 
When $\phi$ is a homomorphism between $C^*$-algebras, 
$K_0(\phi)$ and $K_1(\phi)$ mean the induced homomorphisms on $K$-groups. 

Let $A$ and $B$ be unital $C^*$-algebras. 
We denote by $\Hom(A,B)$ 
the set of all unital homomorphisms from $A$ to $B$. 
Two unital homomorphisms $\phi,\psi\in\Hom(A,B)$ are said to be 
asymptotically unitarily equivalent, 
if there exists a continuous family of unitaries 
$(u_t)_{t\in[0,\infty)}$ in $B$ such that 
\[
\phi(a)=\lim_{t\to\infty}\Ad u_t(\psi(a))
\]
for all $a\in A$. 
When there exists a sequence of unitaries $(u_n)_{n\in\N}$ in $B$ 
such that 
\[
\phi(a)=\lim_{n\to\infty}\Ad u_n(\psi(a))
\]
for all $a\in A$, 
$\phi$ and $\psi$ are said to be approximately unitarily equivalent. 
An automorphism $\alpha\in\Aut(A)$ is said to be 
asymptotically (resp. approximately) inner 
if $\alpha$ is asymptotically (resp. approximately) 
unitarily equivalent to the identity map. 

Let $\alpha:\Gamma\curvearrowright A$ be 
an action of a discrete group $\Gamma$ on a unital $C^*$-algebra $A$. 
The fixed point subalgebra of $A$ is $A^\alpha$. 
The canonical implementing unitaries 
in the reduced crossed product $C^*$-algebra $A\rtimes_\alpha\Gamma$ 
are written by $(\lambda^\alpha_g)_{g\in\Gamma}$. 
The set of all automorphisms $\phi\in\Aut(A\rtimes_\alpha\Gamma)$ satisfying 
$\phi(\lambda^\alpha_g)\lambda^{\alpha*}_g\in A$ for any $g\in\Gamma$ 
is denoted by $\Aut_{\hat\Gamma}(A\rtimes_\alpha\Gamma)$. 
Two automorphisms $\phi,\psi\in\Aut_{\hat\Gamma}(A\rtimes_\alpha\Gamma)$ 
are said to be $\hat\Gamma$-asymptotically unitarily equivalent, 
if there exists a continuous family of unitaries 
$(u_t)_{t\in[0,\infty)}$ in $A$ such that 
\[
\phi(x)=\lim_{t\to\infty}\Ad u_t(\psi(x))
\]
for all $x\in A\rtimes_\alpha\Gamma$. 
In an analogous way, 
one can define $\hat\Gamma$-approximately unitarily equivalence. 
An automorphism $\phi\in\Aut_{\hat\Gamma}(A\rtimes_\alpha\Gamma)$ is said 
to be $\hat\Gamma$-asymptotically (resp. $\hat\Gamma$-approximately) inner 
if $\alpha$ is $\hat\Gamma$-asymptotically (resp. $\hat\Gamma$-approximately) 
unitarily equivalent to the identity map. 

We recall several definitions from \cite{IM,MS}. 

\begin{df}[{\cite[Definition 2.2]{IM}}]\label{asymp}
Let $\Gamma$ be a countable discrete group and 
let $A$ be a unital $C^*$-algebra. 
An action $\alpha:\Gamma\curvearrowright A$ is said to be 
asymptotically representable, 
if there exists a continuous family of unitaries 
$(v_g(t))_{t\in[0,\infty)}$ in $U(A)$ for each $g\in\Gamma$ 
such that 
\[
\lim_{t\to\infty}\lVert v_g(t)v_h(t)-v_{gh}(t)\rVert=0, 
\]
\[
\lim_{t\to\infty}\lVert\alpha_g(v_h(t))-v_{ghg^{-1}}(t)\rVert=0, 
\]
and 
\[
\lim_{t\to\infty}\lVert v_g(t)av_g(t)^*-\alpha_g(a)\rVert=0
\]
hold for all $g,h\in\Gamma$ and $a\in A$. 

Approximate representability is defined in an analogous way 
(see \cite[Definition 3.6]{I04Duke}). 
\end{df}

\begin{df}[{\cite[Definition 2.7]{MS}}]
Let $\alpha:\Gamma\curvearrowright A$ be an action 
of a countable discrete group $\Gamma$ on a unital $C^*$-algebra $A$ 
such that $\tau\circ\alpha_g=\tau$ for any $\tau\in T(A)$ and $g\in\Gamma$. 
We say that $\alpha$ is strongly outer 
if the weak extension of $\alpha_g$ to an automorphism of $\pi_\tau(A)''$ 
is outer for any $g\in\Gamma\setminus\{e\}$ and $\tau\in T(A)$, 
where $\pi_\tau$ is the GNS representation of $A$ 
associated with $\tau$. 
Strong outerness of cocycle actions is defined in the same way. 
\end{df}

\begin{df}[{\cite[Definition 2.1]{MS}}]
Let $\alpha:\Gamma\curvearrowright A$ and $\beta:\Gamma\curvearrowright B$ 
be actions of a countable discrete group $\Gamma$ 
on unital $C^*$-algebras $A$ and $B$. 
\begin{enumerate}
\item The two actions $\alpha$ and $\beta$ are said to be conjugate, 
when there exists an isomorphism $\mu:A\to B$ such that 
$\alpha_g=\mu^{-1}\circ\beta_g\circ\mu$ for all $g\in\Gamma$. 
\item A family of unitaries $(u_g)_{g\in\Gamma}$ in $A$ is called 
an $\alpha$-cocycle, 
if one has $u_g\alpha_g(u_h)=u_{gh}$ for all $g,h\in\Gamma$. 
When $(u_g)_g$ is an $\alpha$-cocycle, 
the perturbed action $\alpha^u:\Gamma\curvearrowright A$ is 
defined by $\alpha^u_g=\Ad u_g\circ\alpha_g$. 
\item The two actions $\alpha$ and $\beta$ are said to be cocycle conjugate, 
if there exists an $\alpha$-cocycle $(u_g)_{g\in\Gamma}$ in $A$ such that 
$\alpha^u$ is conjugate to $\beta$. 
\item The two actions $\alpha$ and $\beta$ are said to be 
strongly cocycle conjugate, 
if there exist an $\alpha$-cocycle $(u_g)_{g\in\Gamma}$ in $A$ 
and a sequence of unitaries $(v_n)_{n=1}^\infty$ in $A$ such that 
$\alpha^u$ is conjugate to $\beta$ and 
$\lim_{n\to\infty}\lVert u_g-v_n\alpha_g(v_n^*)\rVert=0$ for all $g\in\Gamma$. 
\end{enumerate}
\end{df}

Let $A$ be a separable $C^*$-algebra. 
Set 
\[
c_0(A)=\{(a_n)_n\in\ell^\infty(\N,A)\mid
\lim_{n\to\infty}\lVert a_n\rVert=0\},\quad 
A^\infty=\ell^\infty(\N,A)/c_0(A). 
\]
We identify $A$ with the $C^*$-subalgebra of $A^\infty$ 
consisting of equivalence classes of constant sequences. 
We let 
\[
A_\infty=A^\infty\cap A'
\]
and call it the central sequence algebra of $A$. 
A sequence $(x_n)_n\in\ell^\infty(\N,A)$ is called a central sequence 
if $\lVert[a,x_n]\rVert\to0$ as $n\to\infty$ for all $a\in A$. 
A central sequence is a representative of an element in $A_\infty$. 
When $\alpha$ is an automorphism on $A$ or 
an action of a discrete group on $A$, 
we can consider its natural extension on $A^\infty$ and $A_\infty$. 
We denote it by the same symbol $\alpha$. 

We denote the canonical basis of $\Z^N$ by $\xi_1,\xi_2,\dots,\xi_N$, 
that is, 
\[
\xi_i=(0,0,\dots,1,\dots,0,0), 
\]
where $1$ is in the $i$-th component. 
We regard $\Z^{N-1}$ as a subgroup of $\Z^N$ 
via the map $(n_1,n_2,\dots,n_{N-1})\mapsto(n_1,n_2,\dots,n_{N-1},0)$. 

We would like to recall 
the definition of the Rohlin property for $\Z^N$-actions 
on unital $C^*$-algebras (see \cite[Section 2]{N1}). 
Let $\xi_1,\xi_2,\dots,\xi_N$ be the canonical basis of $\Z^N$ as above. 
For $m=(m_1,m_2,\dots,m_N)$ and $n=(n_1,n_2,\dots,n_N)$ in $\Z^N$, 
$m\leq n$ means $m_i\leq n_i$ for all $i=1,2,\dots,N$. 
For $m=(m_1,m_2,\dots,m_N)\in\N^N$, we let 
\[
m\Z^N=\{(m_1n_1,m_2n_2,\dots,m_Nn_N)\in\Z^N
\mid (n_1,n_2,\dots,n_N)\in\Z^N\}. 
\]
For simplicity, we denote $\Z^N/m\Z^N$ by $\Z_m$. 
Moreover, we may identify $\Z_m=\Z^N/m\Z^N$ with 
\[
\{(n_1,n_2,\dots,n_N)\in\Z^N
\mid 0\leq n_i\leq m_i{-}1\quad\forall i=1,2,\dots,N\}. 
\]

\begin{df}\label{Rohlin}
Let $\alpha$ be an action of $\Z^N$ on a unital $C^*$-algebra $A$. 
Then $\alpha$ is said to have the Rohlin property, 
if for any $m\in\N$ there exist $R\in\N$ and 
$m^{(1)},m^{(2)},\dots,m^{(R)}\in\N^N$ 
with $m^{(1)},\dots,m^{(R)}\geq(m,m,\dots,m)$ 
satisfying the following: 
For any finite subset $F$ of $A$ and $\ep>0$, 
there exists a family of projections 
\[
e^{(r)}_g\quad (r=1,2,\dots,R,\quad g\in\Z_{m^{(r)}})
\]
in $A$ such that 
\[
\sum_{r=1}^R\sum_{g\in\Z_{m^{(r)}}}e^{(r)}_g=1,\quad 
\lVert[a,e^{(r)}_g]\rVert<\ep,\quad 
\lVert\alpha_{\xi_i}(e^{(r)}_g)-e^{(r)}_{g+\xi_i}\rVert<\ep
\]
for any $a\in F$, $r=1,2,\dots,R$, $i=1,2,\dots,N$ and $g\in\Z_{m^{(r)}}$, 
where $g+\xi_i$ is understood modulo $m^{(r)}\Z^N$. 
\end{df}

\begin{rem}\label{restate}
Clearly, we can restate the definition of the Rohlin property 
as follows. 
For any $m\in\N$ there exist $R\in\N$, 
$m^{(1)},m^{(2)},\dots,m^{(R)}\in\N^N$ 
with $m^{(1)},\dots,m^{(R)}\geq(m,m,\dots,m)$ and 
a family of projections 
\[
e^{(r)}_g\quad (r=1,2,\dots,R,\quad g\in\Z_{m^{(r)}})
\]
in $A_\infty=A^\infty\cap A'$ such that 
\[
\sum_{r=1}^R\sum_{g\in\Z_{m^{(r)}}}e^{(r)}_g=1,\quad 
\alpha_{\xi_i}(e^{(r)}_g)=e^{(r)}_{g+\xi_i}
\]
for any $r=1,2,\dots,R$, $i=1,2,\dots,N$ and $g\in\Z_{m^{(r)}}$, 
where $g+\xi_i$ is understood modulo $m^{(r)}\Z^N$. 

We can also restate the Rohlin property as follows (\cite[Remark 2]{N1}). 
For any $n,m\in\N$ with $1\leq n\leq N$, 
there exist $R\in\N$, 
natural numbers $m^{(1)},m^{(2)},\dots,m^{(R)}\geq m$ and 
a family of projections 
\[
e^{(r)}_j\quad (r=1,2,\dots,R,\quad j=0,1,\dots,m^{(r)}{-}1)
\]
in $A_\infty=A^\infty\cap A'$ such that 
\[
\sum_{r=1}^R\sum_{j=0}^{m^{(r)}-1}e^{(r)}_j=1,\quad 
\alpha_{\xi_n}(e^{(r)}_j)=e^{(r)}_{j+1},\quad 
\alpha_{\xi_i}(e^{(r)}_j)=e^{(r)}_j
\]
for any $r=1,2,\dots,R$, $i=1,2,\dots,N$ with $i\neq n$ and 
$j=0,1,\dots,m^{(r)}{-}1$, 
where the index $j{+}1$ is understood modulo $m^{(r)}$. 
\end{rem}

\section{$\Z^N$-actions with the Rohlin property}

In this section we prove that 
all $\Z^N$-actions on a UHF algebra of infinite type with the Rohlin property 
are cocycle conjugate to each other (Theorem \ref{CCofRohlin1}). 
We also prove that 
all asymptotically representable $\Z^N$-actions with the Rohlin property 
on a unital simple AH algebra with real rank zero and slow dimension growth 
are cocycle conjugate to each other (Theorem \ref{CCofRohlin2}). 

\begin{lem}\label{beingAT}
Let $\alpha:\Z^N\curvearrowright A$ be 
a strongly outer action of $\Z^N$ on a UHF algebra $A$. 
\begin{enumerate}
\item Let $\alpha':\Z^{N-1}\curvearrowright A$ be the $\Z^{N-1}$-action 
generated by the first $N{-}1$ generators of $\alpha$ and 
let $\tilde\alpha_{\xi_N}$ denote the canonical extension 
of the last generator of $\alpha$ to $A\rtimes_{\alpha'}\Z^{N-1}$. 
Then $K_i(\tilde\alpha_{\xi_N})=\id$ for $i=0,1$. 
\item The crossed product $A\rtimes_\alpha\Z^N$ is 
a unital simple AT algebra of real rank zero with a unique trace. 
\item Let $\iota_N:C^*(\Z^N)\to A\rtimes_\alpha\Z^N$ be 
the canonical inclusion. 
Then 
\[
K_i(\iota_N)\otimes\id_\Q:K_i(C^*(\Z^N))\otimes\Q
\to K_i(A\rtimes_\alpha\Z^N)\otimes\Q
\]
is an isomorphism for $i=0,1$. 
\end{enumerate}
\end{lem}
\begin{proof}
The proof is by induction on $N$. 
When $N=1$, (1) and (3) are clear. 
(2) follows from (1) and \cite[Corollary 5.9]{MS} 
(or \cite[Theorem 1.3]{K95crelle}). 

Suppose that the assertions hold for $N{-}1$. 
From the commutative diagram 
\[
\begin{CD}
C^*(\Z^{N-1}) @>\id>> C^*(\Z^{N-1}) \\
@V\iota_{N-1} VV @V\iota_{N-1} VV \\
A\rtimes_{\alpha'}\Z^{N-1} @>>\tilde\alpha_{\xi_N}> 
A\rtimes_{\alpha'}\Z^{N-1}, 
\end{CD}
\]
one obtains 
$K_i(\iota_{N-1})=K_i(\tilde\alpha_{\xi_N})\circ K_i(\iota_{N-1})$. 
By (3) for $N{-}1$, $K_i(\iota_{N-1})\otimes\id_\Q$ is an isomorphism, 
and so $K_i(\tilde\alpha_{\xi_N})\otimes\id_\Q$ is the identity 
on $K_i(A\rtimes_{\alpha'}\Z^{N-1})\otimes\Q$. 
By (2) for $N{-}1$, $K_i(A\rtimes_{\alpha'}\Z^{N-1})$ is torsion free. 
Therefore $K_i(\tilde\alpha_{\xi_N})=\id$ for $i=0,1$. 
Thus (1) for $N$ has been shown. 

(2) for $N$ follows from (2) for $N{-}1$, (1) for $N$ and 
\cite[Corollary 5.9]{MS}. 

Finally we prove (3). 
Thanks to the naturality of the Pimsner-Voiculescu exact sequence and (1), 
we have the following commutative diagram: 
\[
\begin{CD}
0 @>>> K_i(C^*(\Z^{N-1})) @>>> K_i(C^*(\Z^N)) @>>> 
K_{1-i}(C^*(\Z^{N-1})) @>>> 0 \\
@. @VK_i(\iota_{N-1}) VV @VK_i(\iota_N) VV @VK_{1-i}(\iota_{N-1}) VV @. \\
0 @>>> K_i(A\rtimes_{\alpha'}\Z^{N-1}) @>>> 
K_i(A\rtimes_\alpha\Z^N) @>>> K_{1-i}(A\rtimes_{\alpha'}\Z^{N-1}) @>>> 0
\end{CD}
\]
for $i=0,1$, where the horizontal sequences are exact. 
These two sequences are still exact 
when one takes tensor products with $\Q$. 
By (3) for $N{-}1$, $K_i(\iota_{N-1})\otimes\id_\Q$ is an isomorphism. 
It follows that $K_i(\iota_N)\otimes\id_\Q$ is also an isomorphism. 
\end{proof}

\begin{lem}\label{beingasymp}
Let $\alpha:\Z^N\curvearrowright A$ be 
a strongly outer action of $\Z^N$ on a UHF algebra $A$ of infinite type 
and let $D=K_0(A)$. 
\begin{enumerate}
\item Let $\alpha':\Z^{N-1}\curvearrowright A$ be the $\Z^{N-1}$-action 
generated by the first $N{-}1$ generators of $\alpha$ and 
let $\tilde\alpha_{\xi_N}$ denote the canonical extension 
of the last generator of $\alpha$ to $A\rtimes_{\alpha'}\Z^{N-1}$. 
Then $\tilde\alpha_{\xi_N}$ is asymptotically inner. 
\item For each $i=0,1$, 
$K_i(A\rtimes_\alpha\Z^N)$ is isomorphic to $D^{2^{N-1}}$. 
\item Let $\iota_N:C^*(\Z^N)\to A\rtimes_\alpha\Z^N$ be 
the canonical inclusion. 
Then 
\[
K_i(\iota_N)\otimes\id_D:K_i(C^*(\Z^N))\otimes D
\to K_i(A\rtimes_\alpha\Z^N)\otimes D\cong D^{2^{N-1}}
\]
is an isomorphism for $i=0,1$. 
\end{enumerate}
\end{lem}
\begin{proof}
Note that $D\otimes D\cong D$ and $\Ext(D,D)=0$, 
because $A$ is of infinite type. 
The proof is by induction on $N$. 
When $N=1$, the assertions trivially hold. 

Suppose that the lemma is known for $N{-}1$. 
We would like to show that $\tilde\alpha_{\xi_N}$ is asymptotically inner. 
By Lemma \ref{beingAT}, 
$A\rtimes_{\alpha'}\Z^{N-1}$ is a unital simple AT algebra of real rank zero 
with a unique trace $\tau$ and 
$\tilde\alpha_{\xi_N}$ is approximately inner. 
By virtue of \cite[Theorem 3.1]{KK01JOP}, 
it suffices to show that 
the $\OrderExt$ invariant of $\tilde\alpha_{\xi_N}$ is trivial. 
By (2) for $N{-}1$, 
\[
\Ext(K_i(A\rtimes_{\alpha'}\Z^{N-1}),K_{1-i}(A\rtimes_{\alpha'}\Z^{N-1}))
\cong\Ext(D^{2^{N-2}},D^{2^{N-2}})=0
\]
for $i=0,1$. 
Let 
\[
B=\{f:[0,1]\to A\rtimes_{\alpha'}\Z^{N-1}\mid \tilde\alpha_{\xi_N}(f(0))=f(1)\}
\]
be the mapping torus of $\tilde\alpha_{\xi_N}$. 
We have the short exact sequence: 
\[
\begin{CD}
0 @>>> C_0((0,1),A\rtimes_{\alpha'}\Z^{N-1}) @>>> 
B @>\pi>> A\rtimes_{\alpha'}\Z^{N-1} @>>> 0
\end{CD}
\]
Since the $\Ext$ group is trivial, one has $K_i(B)\cong D^{2^{N-1}}$. 
Let $R:K_1(B)\to\R$ be the rotation map defined in \cite[Lemma 2.1]{KK01JOP}. 
For any unitary $u\in C^*(\Z^{N-1})\otimes M_n$, 
the constant function $\iota_{N-1}(u)$ on the closed interval $[0,1]$ 
belongs to $B\otimes M_n$, 
because $\iota_{N-1}(u)$ is fixed by $\tilde\alpha_{\xi_N}$. 
Hence $R([\iota_{N-1}(u)])=0$, 
where $\iota_{N-1}(u)$ is identified with the constant function on $[0,1]$. 
Therefore there exists a homomorphism 
$\rho:\Ima K_1(\iota_{N-1})\to\Ker R$ such that $K_1(\pi)\circ\rho=\id$. 
By (3) for $N{-}1$, 
\[
K_1(\iota_{N-1})\otimes\id_D:K_1(C^*(\Z^{N-1}))\otimes D\to 
K_1(A\rtimes_{\alpha'}\Z^{N-1})\otimes D
\]
is an isomorphism, and 
$K_1(A\rtimes_{\alpha'}\Z^{N-1})\otimes D$ is naturally 
identified with $K_1(A\rtimes_{\alpha'}\Z^{N-1})$. 
It follows that 
$\rho$ extends to a homomorphism 
from $K_1(A\rtimes_{\alpha'}\Z^{N-1})$ to $\Ker R$ 
satisfying $K_1(\pi)\circ\rho=\id$. 
By \cite[Proposition 2.5]{KK01JOP}, we can conclude that 
the $\OrderExt$ invariant of $\tilde\alpha_{\xi_N}$ is trivial, 
thereby completing the proof of (1) for $N$. 

(2) for $N$ readily follows from the Pimsner-Voiculescu exact sequence 
and $\Ext(D,D)=0$. 
(3) can be shown in a similar way to Lemma \ref{beingAT} (3). 
\end{proof}

\begin{lem}\label{Lipschitz}
Let $\Gamma$ be a countable discrete amenable group and 
let $\alpha:\Gamma\curvearrowright A$ be an approximately representable action 
on a unital $C^*$-algebra $A$. 
Suppose that $A\rtimes_\alpha\Gamma$ is a unital simple AH algebra 
with real rank zero and slow dimension growth. 
For any finite subset $F\subset A\rtimes_\alpha\Gamma$ and $\ep>0$, 
there exist a finite subset $G\subset A\rtimes_\alpha\Gamma$ and $\delta>0$ 
satisfying the following. 
If $u:[0,1]\to A\rtimes_\alpha\Gamma$ is a path of unitaries such that 
\[
u(0)\in A,\quad u(1)\in A,\quad \lVert[a,u(t)]\rVert<\delta
\]
for any $a\in G$ and $t\in[0,1]$, then 
there exists a path of unitaries $w:[0,1]\to A$ such that 
\[
\Lip(w)<11\pi,\quad w(0)=u(0),\quad w(1)=u(1), 
\]
and 
\[
\lVert[a,w(t)]\rVert<\ep
\]
for any $a\in F$ and $t\in[0,1]$. 
\end{lem}
\begin{proof}
Since $\alpha$ is approximately representable, 
we can find a family of unitaries $(v_g)_{g\in\Gamma}$ in $A^\infty$ 
such that 
\[
v_gv_h=v_{gh},\quad \alpha_g(v_h)=v_{ghg^{-1}}
\quad\text{and}\quad v_gav_g^*=\alpha_g(a)
\]
for all $g,h\in\Gamma$ and $a\in A$. 
Define a unital homomorphism $\phi:A\rtimes_\alpha\Gamma\to A^\infty$ by 
\[
\phi(a)=a\quad\text{and}\quad \phi(\lambda^\alpha_g)=v_g
\]
for every $a\in A$ and $g\in\Gamma$. 
It is easy to see that 
$\phi(\lambda^\alpha_gx\lambda^{\alpha*}_g)=\alpha_g(\phi(x))$ holds 
for any $g\in\Gamma$ and $x\in A\rtimes_\alpha\Gamma$. 

Suppose that 
we are given a finite subset $F\subset A\rtimes_\alpha\Gamma$ and $\ep>0$. 
Without loss of generality, we may assume that 
$F$ is of the from $F_0\cup\{\lambda^\alpha_g\mid g\in\Gamma_0\}$, 
where $F_0$ is a finite subset of $A$ and 
$\Gamma_0$ is a finite subset of $\Gamma$. 
Applying \cite[Lemma 3.10]{M10CMP} to $F$ and $\ep>0$, 
we obtain a finite subset $G\subset A\rtimes_\alpha\Gamma$ and $\delta>0$. 
Let $u:[0,1]\to A\rtimes_\alpha\Gamma$ be a path of unitaries satisfying 
$u(0)\in A$, $u(1)\in A$ and 
\[
\lVert[a,u(t)]\rVert<\delta
\]
for any $a\in G$ and $t\in[0,1]$. 
By \cite[Lemma 3.10]{M10CMP}, 
there exists a path of unitaries $w:[0,1]\to A\rtimes_\alpha\Gamma$ such that 
\[
\Lip(w)<11\pi,\quad w(0)=u(0),\quad w(1)=u(1), 
\]
and 
\[
\lVert[a,w(t)]\rVert<\ep
\]
for any $a\in F$ and $t\in[0,1]$. 
Define $\tilde w:[0,1]\to U(A^\infty)$ by $\tilde w(t)=\phi(w(t))$. 
Then we get 
\[
\Lip(\tilde w)<11\pi,\quad \tilde w(0)=u(0),\quad \tilde w(1)=u(1), 
\]
\[
\lVert[a,\tilde w(t)]\rVert<\ep
\]
for any $a\in F$ and $t\in[0,1]$, 
which completes the proof. 
\end{proof}

\begin{thm}\label{CCofRohlin1}
Let $\alpha$ and $\beta$ be $\Z^N$-actions on a UHF algebra of infinite type 
with the Rohlin property. 
Then $\alpha$ and $\beta$ are cocycle conjugate. 
In particular, they are asymptotically representable. 
\end{thm}
\begin{proof}
The proof is by induction on $N$. 
The case $N{=}1$ was shown in \cite[Theorem 1.3]{K95crelle}. 
Suppose that the theorem is known for $N{-}1$. 
Let $A$ be a UHF algebra of infinite type and 
let $\alpha,\beta$ be $\Z^N$-actions on $A$ with the Rohlin property. 
Let $\alpha'$ and $\beta'$ be the $\Z^{N-1}$-actions 
generated by the first $N{-}1$ generators of $\alpha$ and $\beta$, 
respectively. 
From the induction hypothesis, 
by conjugating $\beta$ if necessary, we may assume that 
there exists an $\alpha'$-cocycle $(u_g)_{g\in\Z^{N-1}}$ in $A$ 
such that $\beta'_g=\Ad u_g\circ\alpha'_g$. 
Moreover, $\alpha'$ and $\beta'$ are asymptotically representable. 
It is easy to check 
\[
\beta_{\xi_N}\circ\alpha'_g
=(\Ad\beta_{\xi_N}(u_g^*)u_g)\circ\alpha'_g\circ\beta_{\xi_N}
\]
for all $g\in\Z^{N-1}$ and 
$(\beta_{\xi_N}(u_g^*)u_g)_g$ is an $\alpha'$-cocycle. 
Let $B_\alpha$ and $B_\beta$ be the crossed product of $A$ 
by the $\Z^{N-1}$-actions $\alpha'$ and $\beta'$, respectively. 
By Lemma \ref{beingAT} (2), 
$B_\alpha$ and $B_\beta$ are unital simple AT algebras of real rank zero. 
One can define the isomorphism $\pi:B_\beta\to B_\alpha$ 
by $\pi(a)=a$ for all $a\in A$ and 
$\pi(\lambda^{\beta'}_g)=u_g\lambda^{\alpha'}_g$ for all $g\in\Z^{N-1}$. 
The automorphisms $\alpha_{\xi_N}$ and $\beta_{\xi_N}$ of $A$ 
extend to automorphisms $\tilde\alpha_{\xi_N}$ and $\tilde\beta_{\xi_N}$ 
of $B_\alpha$ and $B_\beta$, respectively. 

We apply the argument of \cite[Theorem 5.1]{K98JFA} to automorphisms 
$\pi\circ\tilde\beta_{\xi_N}\circ\pi^{-1}$ and $\tilde\alpha_{\xi_N}$ 
of $B_\alpha$. 
By Lemma \ref{beingasymp} (1), 
both $\tilde\alpha_{\xi_N}$ and $\tilde\beta_{\xi_N}$ are 
asymptotically inner. 
Hence $\pi\circ\tilde\beta_{\xi_N}\circ\pi^{-1}$ and $\tilde\alpha_{\xi_N}$ 
are asymptotically unitarily equivalent. 
Since $\alpha'$ is asymptotically representable, 
by \cite[Theorem 4.8]{IM}, we can conclude that 
they are $\T^{N-1}$-asymptotically unitarily equivalent. 
Besides 
these two automorphisms have the Rohlin property as single automorphisms 
and the Rohlin projections can be chosen from $(A_\infty)^{\alpha'}$, 
because $\alpha$ and $\beta$ have the Rohlin property as $\Z^N$-actions 
(see Remark \ref{restate}). 
Therefore, 
by using Lemma \ref{Lipschitz} instead of \cite[Lemma 4.4]{K98JFA}, 
the usual intertwining argument shows the following 
(see \cite[Theorem 4.11]{IM} and \cite[Theorem 6.1]{M10CMP} 
for similar arguments): 
There exist an approximately inner automorphism 
$\mu\in\Aut_{\T^{N-1}}(B_\alpha)$ and a unitary $v\in A$ 
such that 
\begin{equation}
\mu\circ\pi\circ\tilde\beta_{\xi_N}\circ\pi^{-1}\circ\mu^{-1}
=\Ad v\circ\tilde\alpha_{\xi_N}. 
\label{rat}
\end{equation}
By restricting this equality to $A$, we get 
\begin{equation}
(\mu|A)\circ\beta_{\xi_N}\circ(\mu|A)^{-1}
=\Ad v\circ\alpha_{\xi_N}. 
\label{cow}
\end{equation}
For each $g\in\Z^{N-1}$, 
let $w_g\in A$ be the unitary satisfying 
$\mu(\lambda^{\alpha'}_g)=w_g\lambda^{\alpha'}_g$. 
Then $(w_g)_g$ is an $\alpha'$-cocycle and 
\begin{equation}
(\mu|A)\circ\beta_g\circ(\mu|A)^{-1}
=(\mu|A)\circ\Ad u_g\circ\alpha_g\circ(\mu|A)^{-1}
=\Ad\mu(u_g)w_g\circ\alpha_g
\label{tiger}
\end{equation}
holds for every $g\in\Z^{N-1}$. 
It is not so hard to see that $(\mu(u_g)w_g)_g$ is also an $\alpha'$-cocycle. 
From \eqref{rat}, one can see that 
\begin{align*}
(\mu\circ\pi\circ\tilde\beta_{\xi_N}\circ\pi^{-1}\circ\mu^{-1})
(\lambda^{\alpha'}_g)
&=(\mu\circ\pi\circ\tilde\beta_{\xi_N}\circ\pi^{-1})
(\mu^{-1}(w_g^*)\lambda^{\alpha'}_g) \\
&=\mu(\beta_{\xi_N}(\mu^{-1}(w_g^*))
\beta_{\xi_N}(u_g^*)u_g\lambda^{\alpha'}_g) \\
&=(\Ad v\circ\alpha_{\xi_N})(w_g^*\mu(u_g^*))\mu(u_g)w_g\lambda^{\alpha'}_g \\
&=v\alpha_{\xi_N}(w_g^*\mu(u_g^*))v^*\mu(u_g)w_g\lambda^{\alpha'}_g
\end{align*}
is equal to 
\[
(\Ad v\circ\tilde\alpha_{\xi_N})(\lambda^{\alpha'}_g)
=v\lambda^{\alpha'}_gv^*
=v\alpha_g(v^*)\lambda^{\alpha'}_g
\]
for any $g\in\Z^{N-1}$. 
Hence one obtains 
\[
v\alpha_{\xi_N}(\mu(u_g)w_g)=\mu(u_g)w_g\alpha_g(v). 
\]
Thus, $(\mu(u_g)w_g)_g$ and $v$ give rise to an $\alpha$-cocycle. 
It follows from \eqref{cow} and \eqref{tiger} that 
$\alpha$ and $\beta$ are cocycle conjugate. 

There exists an asymptotically representable $\Z^N$-action on $A$ 
with the Rohlin property, 
and so any $\Z^N$-action on $A$ with the Rohlin property is 
asymptotically representable. 
\end{proof}

\begin{lem}
Let $A$ be a unital simple AH algebra 
with real rank zero and slow dimension growth. 
Then there exists an asymptotically inner automorphism $\sigma\in\Aut(A)$ 
with the Rohlin property such that 
the crossed product $A\rtimes_\sigma\Z$ is again 
a unital simple AH algebra with real rank zero and slow dimension growth. 
\end{lem}
\begin{proof}
By \cite{EGL}, 
one can find an increasing sequence $\{A_n\}_n$ of unital subalgebras of $A$ 
such that the following hold. 
\begin{itemize}
\item $\bigcup_nA_n$ is dense in $A$. 
\item $A_n$ is of the form 
$\bigoplus_{i=1}^{k_n}p_{n,i}M_{l(n,i)}(C(X_{n,i}))p_{n,i}$, 
where $X_{n,i}$ is  a compact Hausdorff space 
with topological dimension at most three and $p_{n,i}$ is a projection. 
\item There exists a unital embedding 
$\pi_n:M_n\oplus M_{n+1}\to A_{n+1}\cap A_n'$. 
\end{itemize}
Let $x_n$ be a unitary of $M_n(\C)$ such that 
$\Sp(x_n)=\{\zeta^k\mid k=0,1,\dots,n{-}1\}$, 
where $\zeta=\exp(2\pi\sqrt{-1}/n)$. 
Let $y_n=\pi_n(x_n\oplus x_{n+1})$. 
Define an automorphism $\sigma$ of $A$ 
by $\sigma=\lim_{n\to\infty}\Ad(y_1y_2\dots y_n)$. 
Then $\sigma$ is 
an asymptotically inner automorphism with the Rohlin property. 
It is easy to see that $A\rtimes_\sigma\Z$ is an inductive limit of 
the $C^*$-algebras $A_n\otimes C(\T)$. 
Thus $A\rtimes_\sigma\Z$ is a unital AH algebra with slow dimension growth. 
Simplicity also follows because the action $\sigma$ is outer. 
By \cite[Theorem 4.5]{OP}, it has real rank zero. 
The proof is completed. 
\end{proof}

\begin{lem}\label{beingAH}
Let $A$ be a unital simple AH algebra 
with real rank zero and slow dimension growth. 
Let $\alpha:\Z^N\curvearrowright A$ be 
an asymptotically representable action of $\Z^N$ with the Rohlin property. 
Then the crossed product $A\rtimes_\alpha\Z^N$ is 
a unital simple AH algebra with real rank zero and slow dimension growth. 
\end{lem}
\begin{proof}
The proof is by induction on $N$. 
Let $\alpha\in\Aut(A)$ be 
an asymptotically inner automorphism with the Rohlin property. 
By the lemma above, 
there exists an asymptotically inner automorphism $\sigma\in\Aut(A)$ 
with the Rohlin property such that 
the crossed product $A\rtimes_\sigma\Z$ is again 
a unital simple AH algebra with real rank zero and slow dimension growth. 
By \cite[Theorem 4.9]{M10CMP}, 
(the $\Z$-actions generated by) $\alpha$ and $\sigma$ are cocycle conjugate. 
In particular, $A\rtimes_\alpha\Z$ is isomorphic to $A\rtimes_\sigma\Z$. 

Suppose that the lemma is known for $N{-}1$. 
Let $\alpha:\Z^N\curvearrowright A$ be 
an asymptotically representable action of $\Z^N$ with the Rohlin property. 
Let $\alpha'$ be the $\Z^{N-1}$-action 
generated by the first $N{-}1$ generators of $\alpha$. 
From the induction hypothesis, $A\rtimes_{\alpha'}\Z^{N-1}$ is 
a unital simple AH algebra with real rank zero and slow dimension growth. 
Let $\tilde\alpha_{\xi_N}\in\Aut(A\rtimes_{\alpha'}\Z^{N-1})$ be 
the natural extension of the automorphism $\alpha_{\xi_N}$ of $A$. 
Since $\alpha$ is asymptotically representable and 
has the Rohlin property as a $\Z^N$-action, 
$\tilde\alpha_{\xi_N}$ is asymptotically inner and 
has the Rohlin property as a single automorphism. 
By the same argument as above, 
we can conclude that 
\[
A\rtimes_\alpha\Z^N
\cong(A\rtimes_{\alpha'}\Z^{N-1})\rtimes_{\tilde\alpha_{\xi_N}}\Z
\]
is a unital simple AH algebra with real rank zero and slow dimension growth. 
\end{proof}

\begin{thm}\label{CCofRohlin2}
Let $A$ be a unital simple AH algebra 
with real rank zero and slow dimension growth and 
let $\alpha,\beta$ be asymptotically representable $\Z^N$-actions on $A$ 
with the Rohlin property. 
Then $\alpha$ is cocycle conjugate to $\beta$. 
\end{thm}
\begin{proof}
One can prove this in a similar fashion to Theorem \ref{CCofRohlin1}, 
by using Lemma \ref{beingAH} instead of Lemma \ref{beingAT} (2). 
We omit the proof. 
\end{proof}

\section{Rohlin type theorem}

In this section we prove that 
any strongly outer action of $\Z^N$ on a UHF algebra of infinite type 
has the Rohlin property (Theorem \ref{Rohlintype1}). 
We also prove that 
any strongly outer, approximately representable action of $\Z^N$ 
on a unital simple AH algebra with real rank zero, slow dimension growth and 
finitely many extremal traces has the Rohlin property 
(Theorem \ref{Rohlintype2}). 
Combining this with Theorem \ref{CCofRohlin2}, we can conclude that 
all strongly outer, asymptotically representable actions of $\Z^N$ 
on such a unital simple AH algebra are cocycle conjugate 
(Theorem \ref{CCofasymp}). 

\begin{lem}\label{key}
Let $\Gamma$ be a countable discrete amenable group and 
let $\alpha:\Gamma\curvearrowright A$ be an approximately representable action 
on a unital $C^*$-algebra $A$. 
Suppose that 
$\beta\in\Aut_{\hat\Gamma}(A\rtimes_\alpha\Gamma)$ is approximately inner. 
Then for any separable subset $C\subset A_\infty$, 
there exists a unitary $u\in(A_\infty)^\alpha$ such that 
$\beta(x)=uxu^*$ for all $x\in C$. 
\end{lem}
\begin{proof}
By \cite[Remark 4.9]{IM}, $\beta$ is $\hat\Gamma$-approximately inner. 
Thus there exists a sequence $(v_n)_n$ of unitaries of $A$ such that 
\[
\beta(x)=\lim_{n\to\infty}v_nxv_n^*
\]
holds for any $x\in A\rtimes_\alpha\Gamma$. 
Then one can prove the assertion 
in the same way as \cite[Lemma 4.3]{M10CMP}. 
\end{proof}

\begin{lem}\label{cyclicRohlin}
For any $N\in\N$, $l\in\N$ and $\ep>0$, 
there exist $m\in\N^N$ and $k\in\N$ such that the following hold: 
Let $\alpha:\Z^{N-1}\curvearrowright A$ be 
an approximately representable action of $\Z^{N-1}$ 
on a unital $C^*$-algebra $A$. 
Suppose that 
$\beta\in\Aut_{\T^{N-1}}(A\rtimes_\alpha\Z^{N-1})$ is approximately inner. 
Let $e\in A_\infty$ be a projection satisfying 
\[
e\alpha_g(\beta^j(e))=0\quad\forall (g,j)\in\Z_m\setminus\{(0,0)\}, 
\]
where $\Z_m$ is regarded as a subset of $\Z^N$. 
Then there exists a projection $p\in A_\infty$ satisfying the following. 
\begin{enumerate}
\item $\lVert p-\alpha_{\xi_i}(p)\rVert<\ep$ for any $i=1,2,\dots,N{-}1$. 
\item $p\beta^j(p)=0$ for any $j=1,2,\dots,l{-}1$. 
\item $\lVert p-\beta^l(p)\rVert<\ep$. 
\item $\sum_{j=0}^{l-1}\beta^j(p)\leq\sum_{(g,j)\in\Z_m}\alpha_g(\beta^j(e))$. 
\item $[p]$ is equal to $k[e]$ in $K_0(A_\infty)$ and 
$kl\geq(1{-}\ep)\lvert\Z_m\rvert$. 
\end{enumerate}
\end{lem}
\begin{proof}
The proof is by induction on $N$. 
The case $N{=}1$ was shown in the same way as \cite[Lemma 4.3]{K95crelle} 
(see also \cite[Theorem 4.4]{M10CMP}). 
Suppose that the lemma is known for $N{-}1$. 
Suppose that we are given $l\in\N$ and $\ep>0$. 
Choose $k_1,k_2\in\N$ so that 
\[
\frac{1}{k_1}+\frac{1}{\sqrt{k_1}}<\ep\quad\text{and}\quad 
\frac{l(k_1+k_2)}{l(2k_1+k_2-1)+1}>\sqrt{1-\ep}. 
\]
Choose $\delta>0$ so that 
$(2k_1{+}k_2{-}1)\delta<\ep$ and $\delta<1{-}\sqrt{1-\ep}$. 
Applying the lemma to $N{-}1$, $l{=}1$ and $\delta>0$, 
we get $m'\in\N^{N-1}$ and $k'\in\N$. 
We would like to show that 
$m=(m',l(2k_1{+}k_2{-}1)+1)\in\N^N$ and $k'(k_1+k_2)\in\N$ 
meet the requirement. 

Suppose that we are given an approximately representable action 
$\alpha:\Z^{N-1}\curvearrowright A$ on a unital $C^*$-algebra $A$ 
and an approximately inner automorphism 
$\beta\in\Aut_{\T^{N-1}}(A\rtimes_\alpha\Z^{N-1})$. 
Note that for any $g\in\Z^{N-1}$ and $x\in A_\infty$ one has 
\[
\beta(\alpha_g(x))=\beta(\lambda^\alpha_gx\lambda^{\alpha*}_g)
=\beta(\lambda^\alpha_g)\lambda^{\alpha*}_g\lambda^\alpha_g\beta(x)
\lambda^{\alpha*}_g\lambda^\alpha_g\beta_g(\lambda^{\alpha*}_g)
=\alpha_g(\beta(x)), 
\]
because $\beta(\lambda^\alpha_g)\lambda^{\alpha*}_g$ is in $A$. 
Let $e\in A_\infty$ be a projection satisfying 
\[
e\alpha_g(\beta^j(e))=0\quad\forall (g,j)\in\Z_m\setminus\{(0,0)\}, 
\]
where $\Z_m$ is regarded as a subset of $\Z^N$. 
Notice that 
the $\Z^{N-2}$-action generated by the first $N{-}2$ generators of $\alpha$ 
is approximately representable and that 
the canonical extension of $\alpha_{\xi_{N-1}}$ 
to the crossed product by the $\Z^{N-2}$-action is approximately inner. 
By the induction hypothesis, 
we can find a projection $q\in A_\infty$ satisfying the following. 
\begin{itemize}
\item $\lVert q-\alpha_{\xi_i}(q)\rVert<\delta$ for any $i=1,2,\dots,N{-}1$. 
\item $q\leq\sum_{g\in\Z_{m'}}\alpha_g(e)$. 
\item $[q]$ is equal to $k'[e]$ in $K_0(A_\infty)$ and 
$k'\geq(1{-}\delta)\lvert\Z_{m'}\rvert$. 
\end{itemize}
We remark that the second condition implies 
\[
q\beta^j(q)=0\quad\forall j=1,2,\dots,l(2k_1{+}k_2{-}1). 
\]
We construct the desired projection $p$ 
in the same way as \cite[Lemma 2.1]{K95crelle}. 
By using Lemma \ref{key} for $C=\{q\}$, 
we obtain a unitary $u\in(A_\infty)^\alpha$ such that $\beta(q)=uqu^*$. 
For $s,t\in\Z$ with $0\leq s<t$, we set 
\[
u_{s,t}=\beta^{t-1}(u)\dots\beta^{s+2}(u)\beta^{s+1}(u)\beta^s(u)\beta^s(q)
\in A_\infty. 
\]
Then $u_{s,t}$ satisfies 
$u_{s,t}^*u_{s.t}=\beta^s(q)$, $u_{s,t}u_{s,t}^*=\beta^t(q)$ 
and $\lVert u_{s,t}-\alpha_{\xi_i}(u_{s,t})\rVert<\delta$ 
for any $i=1,2,\dots,N{-}1$. 
Define 
\begin{align*}
p=& \sum_{i=1}^{k_1-1}
\left(\frac{i}{k_1}\beta^{l(i-1)}(q)
+\frac{k_1-i}{k_1}\beta^{l(k_1+k_2+i-1)}(q)\right. \\
&\left.+\frac{\sqrt{i(k_1-i)}}{k_1}
\left(u_{l(i-1),l(k_1+k_2+i-1)}+u_{l(i-1),l(k_1+k_2+i-1)}^*\right)\right) \\
&+\sum_{i=k_1}^{k_1+k_2}\beta^{l(i-1)}(q). 
\end{align*}
It is easy to see 
\[
\lVert p-\alpha_{\xi_i}(p)\rVert<(2(k_1-1)+k_2+1)\delta<\ep
\]
for any $i=1,2,\dots,N{-}1$. 
From the construction, 
clearly we have $p\beta^j(p)=0$ for any $j=1,2,\dots,l{-}1$ and 
\[
\lVert p-\beta^l(p)\rVert<\frac{1}{k_1}+\frac{1}{\sqrt{k_1}}<\ep. 
\]
Furthermore, 
\[
\sum_{j=0}^{l-1}\beta^j(p)
\leq\sum_{i=0}^{l(2k_1+k_2-1)}\beta^i(q)
\leq\sum_{(g,j)\in\Z_m}\alpha_g(\beta^j(e)). 
\]
Finally, 
\[
[p]=(k_1+k_2)[q]=k'(k_1+k_2)[e]=k[e]
\]
in $K_0(A_\infty)$ and 
\[
kl=k'(k_1+k_2)l>(1{-}\delta)\lvert\Z_{m'}\rvert\cdot(1{-}\ep)^{1/2}
(l(2k_1{+}k_2{-}1)+1)>(1{-}\ep)\lvert\Z_m\rvert. 
\]
\end{proof}

\begin{lem}\label{tracialRohlin}
Let $A$ be a unital simple separable $C^*$-algebra with tracial rank zero 
and suppose that $A$ has finitely many extremal tracial states. 
Let $(\alpha,u)$ be a strongly outer cocycle action of $\Z^N$ on $A$ 
such that $\tau\circ\alpha_g=\tau$ for any $\tau\in T(A)$ and $g\in\Z^N$. 
Then, for any $m\in\N^N$, 
there exists a central sequence of projections $(e_n)_n$ in $A$ such that 
\[
\lim_{n\to\infty}\tau(e_n)=\lvert\Z_m\rvert^{-1}
\]
for all $\tau\in T(A)$ and 
\[
\lim_{n\to\infty}\lVert\alpha_g(e_n)\alpha_h(e_n)\rVert=0
\]
for all $g\neq h$ in $\Z_m$. 
\end{lem}
\begin{proof}
One can prove this in a similar fashion to \cite[Theorem 3.4]{MS}, 
using \cite[Proposition 4.1]{M10CMP} instead of \cite[Proposition 3.3]{MS}. 
We omit the detail. 
\end{proof}

\begin{lem}\label{equivRohlin}
Let $A$ be a unital simple AH algebra 
with real rank zero and slow dimension growth 
and suppose that $A$ has finitely many extremal tracial states. 
Let $\alpha:\Z^{N-1}\curvearrowright A$ be 
an approximately representable action of $\Z^{N-1}$ and 
let $\beta\in\Aut_{\T^{N-1}}(A\rtimes_\alpha\Z^{N-1})$ be 
an approximately inner automorphism. 
Suppose that the cocycle $\Z^N$-action 
generated by $\alpha$ and $\beta$ on $A$ is strongly outer. 
Then for any $m\in\N$, 
there exist projections $e$ and $f$ in $(A_\infty)^\alpha$ such that 
\[
\beta^m(e)=e,\quad \beta^{m+1}(f)=f
\]
and 
\[
\sum_{j=0}^{m-1}\beta^j(e)+\sum_{j=0}^m\beta^j(f)=1. 
\]
\end{lem}
\begin{proof}
We first prove the following claim. 
For any $m\in\N$, there exist projections $p,q\in(A_\infty)^\alpha$ and 
a partial isometry $w\in(A_\infty)^\alpha$ such that 
\[
w^*w=q,\quad ww^*\leq p,\quad q+\sum_{j=0}^{m-1}\beta^j(p)=1
\quad\text{and}\quad \beta^m(p)=p. 
\]
By virtue of Lemma \ref{cyclicRohlin} and Lemma \ref{tracialRohlin}, 
we can find projections $p,q\in(A_\infty)^\alpha$ such that 
\[
q+\sum_{j=0}^{m-1}\beta^j(p)=1,\quad \beta^m(p)=p
\]
and 
\[
\lim_{n\to\infty}\tau(p_n)=1/m
\]
for any $\tau\in T(A)$, 
where $(p_n)_n$ is a central sequence of projections representing $p$. 
For each $i=1,2,\dots,N{-}1$, 
there exists a sequence of unitaries $(u_{i,n})_n$ in $A$ such that 
$u_{i,n}\to1$ as $n\to\infty$ and 
$u_{i,n}\alpha_{\xi_i}(p_n)u_{i,n}^*=p_n$ for any $n$. 
The cocycle $\Z^{N-1}$-action on $p_nAp_n$ 
generated by $\Ad u_{i,n}\circ\alpha_{\xi_i}$ is strongly outer. 
Let $k\in\N$ and set $K=\Z_{(k,k,\dots,k)}\subset\Z^{N-1}$. 
One can apply Lemma \ref{tracialRohlin} and 
obtain a central sequence $(\tilde p_n)_n$ of projections in $A$ such that 
\[
\tilde p_n\leq p_n,\quad 
\lim_{n\to\infty}\tau(\tilde p_n)=\frac{1}{mk^{N-1}}
\quad\text{and}\quad 
\lim_{n\to\infty}\lVert\tilde p_n\alpha_g(\tilde p_n)\rVert=0
\]
for all $\tau\in T(A)$ and $g\in K\setminus\{0\}$. 
Let $\tilde p\in A_\infty$ be the image of $(\tilde p_n)_n$. 
By \cite[Lemma 3.3]{M10CMP}, 
there exists a partial isometry $v\in A_\infty$ such that 
$v^*v=q$ and $vv^*\leq\tilde p$. 
We define a partial isometry $\tilde v\in A_\infty$ by 
\[
w=\frac{1}{\sqrt{k^{N-1}}}\sum_{g\in K}\alpha_g(v). 
\]
Then one has 
\[
w^*w=q,\quad ww^*\leq p\quad\text{and}\quad 
\lVert w-\alpha_{\xi_i}(w)\rVert<2/\sqrt{k}
\]
for any $i=1,2,\dots,N{-}1$. 
By a standard trick on central sequences, 
we may assume that $w$ belongs to $(A_\infty)^\alpha$, 
thereby completing the proof of the claim. 

We prove the lemma. 
Suppose that we are given $m\in\N$. 
Let $k,l$ be sufficiently large natural numbers. 
By the claim above, 
we can find projections $p,q\in(A_\infty)^\alpha$ and 
a partial isometry $w\in(A_\infty)^\alpha$ such that 
\[
w^*w=q,\quad ww^*\leq p,\quad q+\sum_{j=0}^{klm-1}\beta^j(p)=1
\quad\text{and}\quad \beta^{klm}(p)=p. 
\]
Define $\tilde p,\tilde w\in(A_\infty)^\alpha$ by 
\[
\tilde p=\sum_{j=0}^{k-1}\beta^{jlm}(p)\quad\text{and}\quad 
\tilde w=\frac{1}{\sqrt{k}}\sum_{j=0}^{k-1}\beta^{jlm}(w). 
\]
Then $\tilde p$ is a projection and 
$\tilde w$ is a partial isometry satisfying 
\[
q+\sum_{j=0}^{lm-1}\beta^j(\tilde p)=1,\quad 
\beta^{lm}(\tilde p)=\tilde p
\]
and 
\[
\tilde w^*\tilde w=q,\quad \tilde w\tilde w^*\leq \tilde p,\quad 
\lVert\beta^{lm}(\tilde w)-\tilde w\rVert\leq\frac{2}{\sqrt{k}}. 
\]
Let $D$ be the $C^*$-algebra generated 
by $\tilde w,\beta(\tilde w),\dots,\beta^{lm-1}(\tilde w)$. 
Then $D$ is isomorphic to $M_{lm+1}$ and 
the unit $1_D$ of $D$ is equal to 
$q+\tilde w\tilde w^*+\dots+\beta^{lm-1}(\tilde w\tilde w^*)$. 
From the spectral property of $\beta$ restricted to $D$, 
if $k$ and $l$ are sufficiently large, 
we can obtain projections $e_0,\dots,e_{m-1},f_0,\dots,f_m$ of $D$ such that 
\[
\sum_{i=1}^{m-1}e_i+\sum_{i=1}^mf_i=1_D,\quad 
\beta(e_i)\approx e_{i+1},\quad \beta(f_i)\approx f_{i+1},
\]
where $e_m=e_0$ and $f_{m+1}=f_0$ 
(see \cite{K95crelle,K96JFA} for details). 
We define projections $e'_i$ in $(A_\infty)^\alpha$ by 
\[
e'_i=e_i+\sum_{j=0}^{l-1}\beta^{i+jm}(\tilde p-\tilde w\tilde w^*). 
\]
Then the projections $e'_0,\dots,e'_{m-1},f_0,\dots,f_m\in(A_\infty)^\alpha$ 
meet the requirement approximately. 
The usual reindexation trick completes the proof. 
\end{proof}

\begin{thm}\label{Rohlintype1}
Let $A$ be a UHF algebra of infinite type and 
let $\alpha:\Z^N\curvearrowright A$ be an action of $\Z^N$. 
The following are equivalent. 
\begin{enumerate}
\item $\alpha$ has the Rohlin property. 
\item $\alpha$ is uniformly outer. 
\item $\alpha$ is strongly outer. 
\end{enumerate}
\end{thm}
\begin{proof}
(1)$\Rightarrow$(2) is obvious. 
(2)$\Rightarrow$(3) follows from \cite[Lemma 4.4]{K96JFA}. 
We prove (3)$\Rightarrow$(1). 
The proof is by induction on $N$. 
The case $N{=}1$ was shown in \cite{K95crelle}. 
Suppose that the assertion is known for $N{-}1$. 
Let $\alpha:\Z^N\curvearrowright A$ be a strongly outer action of $\Z^N$. 
Let $\alpha'$ be the $\Z^{N-1}$-action 
generated by the first $N{-}1$ generators of $\alpha$. 
From the induction hypothesis, $\alpha'$ has the Rohlin property. 
It follows from Theorem \ref{CCofRohlin1} that 
$\alpha'$ is asymptotically representable. 
Let $\tilde\alpha_{\xi_N}\in\Aut_{\T^{N-1}}(A\rtimes_{\alpha'}\Z^{N-1})$ 
be the natural extension of the automorphism $\alpha_{\xi_N}$ of $A$. 
By Lemma \ref{beingAT} (1) (or Lemma \ref{beingasymp} (1)), 
$\tilde\alpha_{\xi_N}$ is approximately inner. 
Hence we can apply Lemma \ref{equivRohlin} 
to $\alpha'$ and $\tilde\alpha_{\xi_N}$ and 
obtain Rohlin projections for $\alpha_{\xi_N}$ in $(A_\infty)^{\alpha'}$. 
The same argument works for other generators $\xi_i$ instead of $\xi_N$. 
By Remark \ref{restate}, we can conclude that 
$\alpha$ has the Rohlin property. 
\end{proof}

\begin{rem}
In the proof above 
we have shown the uniqueness of $\alpha$ up to cocycle conjugacy. 
In the next section 
it will be strengthened to strong cocycle conjugacy (Theorem \ref{SCConUHF}). 
\end{rem}

\begin{thm}\label{Rohlintype2}
Let $A$ be a unital simple AH algebra 
with real rank zero and slow dimension growth 
and suppose that $A$ has finitely many extremal tracial states. 
Let $\alpha:\Z^N\curvearrowright A$ be 
an approximately representable action of $\Z^N$. 
The following are equivalent. 
\begin{enumerate}
\item $\alpha$ has the Rohlin property. 
\item $\alpha$ is uniformly outer. 
\item $\alpha$ is strongly outer. 
\end{enumerate}
\end{thm}
\begin{proof}
One can prove this in the same way as Theorem \ref{Rohlintype1}. 
\end{proof}

The following is an immediate consequence of 
Theorem \ref{CCofRohlin2} and Theorem \ref{Rohlintype2}. 

\begin{thm}\label{CCofasymp}
Let $A$ be a unital simple AH algebra 
with real rank zero and slow dimension growth 
and suppose that $A$ has finitely many extremal tracial states. 
Let $\alpha,\beta:\Z^N\curvearrowright A$ be 
strongly outer, asymptotically representable actions of $\Z^N$. 
Then they are cocycle conjugate. 
\end{thm}

\section{Cohomology vanishing}

In this section, we prove a cohomology vanishing theorem 
(Lemma \ref{admissible}). 
As an application, 
we show Theorem \ref{SCConUHF} and Theorem \ref{SCCofasymp}. 

For a unital $C^*$-algebra $A$, 
we say that $K_0(A)$ has no infinitesimal 
if for any $x\in K_0(A)\setminus\{0\}$ 
there exists $\tau\in T(A)$ such that $\tau(x)\neq0$. 

In the next lemma, we use $KK$-theory. 
Let $A$, $B$ and $C$ be $C^*$-algebras. 
For a homomorphism $\phi:A\to B$, 
$KK(\phi)$ means the induced element in $KK(A,B)$. 
We write $KK(\id_A)=1_A$. 
For $x\in KK(A,B)$ and $i=0,1$, 
we let $K_i(x)$ denote the homomorphism from $K_i(A)$ to $K_i(B)$ 
induced by $x$. 
For $x\in KK(A,B)$ and $y\in KK(B,C)$, 
we denote the Kasparov product by $x\cdot y\in KK(A,C)$. 

Let $\Gamma$ be a countable discrete amenable group. 
We denote the canonical generators in $C^*(\Gamma)$ 
by $(\lambda_g)_{g\in\Gamma}$. 
Let $\alpha:\Gamma\curvearrowright A$ be an action of $\Gamma$ 
on a unital $C^*$-algebra $A$. 
We let $\Hom_{\hat\Gamma}(C^*(\Gamma),A\rtimes_{\alpha}\Gamma)$ denote 
the set of all $\phi\in\Hom(C^*(\Gamma),A\rtimes_{\alpha}\Gamma)$ 
such that $\phi(\lambda_g)\lambda^{\alpha*}_g\in A$ for any $g\in\Gamma$. 

\begin{lem}\label{KKcomputation}
Let $\gamma:\Z^N\curvearrowright\mathcal{Z}$ be 
an action of $\Z^N$ on the Jiang-Su algebra and 
let $A$ be a unital $C^*$-algebra such that 
$K_0(A)$ has no infinitesimal and $K_1(A)=0$. 
Then for any $\phi,\psi\in
\Hom_{\T^N}(C^*(\Z^N),(A\otimes\mathcal{Z})\rtimes_{\id\otimes\gamma}\Z^N)$, 
one has $KK(\phi)=KK(\psi)$. 
\end{lem}
\begin{proof}
Let $\iota:C^*(\Z^N)\to\mathcal{Z}\rtimes_\gamma\Z^N$ be 
the canonical embedding. 
In the same way as \cite[Lemma 5.3]{IM}, 
one can see that $KK(\iota)$ gives a $KK$-equivalence. 
It follows that $KK(\id_A\otimes\iota)$ also gives a $KK$-equivalence 
between $A\otimes C^*(\Z^N)$ and 
$A\otimes(\mathcal{Z}\rtimes_\gamma\Z^N)
\cong(A\otimes\mathcal{Z})\rtimes_{\id\otimes\gamma}\Z^N$. 

Let $S=C_0((0,1))$ and let $\mathcal{T}=\C\oplus S$. 
Put $[N]=\{1,2,\dots,N\}$. 
The $N$-fold tensor product $\mathcal{T}^{\otimes N}$ 
has $2^N$ direct sum components and 
each of them is isomorphic to a tensor product of several copies of $S$. 
For $I\subset[N]$, 
we let $S_I\subset\mathcal{T}^{\otimes N}$ denote 
the tensor product of $S$'s 
of the $i$-th tensor product component for all $i\in I$, 
so that 
\[
\mathcal{T}^{\otimes N}=\bigoplus_{I\subset[N]}S_I. 
\]
Note that $S_\emptyset$ is isomorphic to $\C$. 
Let $z\in KK(C^*(\Z),\mathcal{T})$ be an invertible element. 
We denote the $N$-fold tensor product of $z$ 
by $z_N\in KK(C^*(\Z^N),\mathcal{T}^{\otimes N})$. 

Take $\phi\in
\Hom_{\T^N}(C^*(\Z^N),(A\otimes\mathcal{Z})\rtimes_{\id\otimes\gamma}\Z^N)$. 
Define $a\in KK(\mathcal{T}^{\otimes N},A\otimes\mathcal{T}^{\otimes N})$ by 
\[
a=z_N^{-1}\cdot KK(\phi)\cdot(1_A\otimes KK(\iota)^{-1})
\cdot(1_A\otimes z_N). 
\]
As in \cite[Section 6]{IM}, under the identification 
\[
KK(\mathcal{T}^{\otimes N},A\otimes\mathcal{T}^{\otimes N})
=\bigoplus_{I,J\subset[N]}KK(S_I,A\otimes S_J), 
\]
we denote the $KK(S_I,A\otimes S_J)$ component of $a$ by $a(I,J)$. 
By \cite[Lemma 6.15, 6.16]{IM}, 
for each $K\subset[N]$, there exists $b_K\in KK(S_K,A)$ such that 
\[
a(I,J)=\begin{cases}b_{I\setminus J}\otimes1_{S_J} & J\subset I \\
0 & \text{otherwise}\end{cases}
\]
and $b_\emptyset=KK(h)$, where $h:\C\to A$ is the unital homomorphism. 
We would like to show $b_K=0$ for $K\neq\emptyset$. 
When $\lvert K\rvert$ is odd, $KK(S_K,A)\cong K_1(A)$ is zero, 
and so $b_K$ is zero. 
Suppose that $\lvert K\rvert$ is even (and nonzero). 
For any tracial state $\tau$ on 
$(A\otimes\mathcal{Z})\rtimes_{\id\otimes\gamma}\Z^N$ 
and $x\in K_0(A\otimes\mathcal{T}^{\otimes N})$, 
it is easy to see that 
\[
(\tau\circ K_0((1_A\otimes z_N^{-1})\cdot(1_A\otimes KK(\iota))))(x)
=\tau(x_\emptyset), 
\]
where $x_\emptyset\in K_0(A)$ is the $K_0(A\otimes S_\emptyset)$ summand of 
$x\in K_0(A\otimes\mathcal{T}^{\otimes N})$ and 
$A\otimes S_\emptyset$ is identified with $A$. 
Therefore, letting $\omega_K$ be the generator of $K_0(S_K)$, we have 
\begin{align*}
&(\tau\circ K_0(a\cdot(1_A\otimes z_N^{-1})\cdot(1_A\otimes KK(\iota))))
(\omega_K) \\
&=\sum_{J\subset K}
(\tau\circ K_0((1_A\otimes z_N^{-1})\cdot(1_A\otimes KK(\iota)))
\circ K_0(a(K,J)))(\omega_K) \\
&=(\tau\circ K_0((1_A\otimes z_N^{-1})\cdot(1_A\otimes KK(\iota)))
\circ K_0(a(K,\emptyset)))(\omega_K) \\
&=(\tau\circ K_0(b_K))(\omega_K), 
\end{align*}
where the restriction of $\tau$ to $A$ is also denoted by $\tau$. 
On the other hand, $\tau\circ\phi$ is a trace on $C^*(\Z^N)\cong C(\T^N)$, 
and so 
\[
(\tau\circ K_0(\phi)\circ K_0(z_N^{-1}))(\omega_K)=0. 
\]
Combining these equalities, we get $(\tau\circ K_0(b_K))(\omega_K)=0$ 
for any $\tau$. 
Since $K_0(A)$ has no infinitesimal, $K_0(b_K)(\omega_K)$ is zero. 
It follows that $K_0(b_K)$ is zero, 
because $\omega_K$ is the generator of $K_0(S_K)$. 
Hence $b_K$ is zero. 
We have thus shown the lemma. 
\end{proof}

In the next lemma, we let $\mathcal{C}_0$ denote 
the class of unital simple stably finite $C^*$-algebras 
introduced in \cite[Definition 2.5]{MS}. 

\begin{lem}\label{Z}
Let $\gamma:\Z^N\curvearrowright\mathcal{Z}$ be 
a strongly outer action of $\Z^N$ on the Jiang-Su algebra. 
Then $\mathcal{Z}\rtimes_\gamma\Z^N$ belongs to $\mathcal{C}_0$ 
and has a unique tracial state. 
\end{lem}
\begin{proof}
The proof is by induction on $N$. 
Assume that the lemma is known for $N{-}1$. 
Let $\gamma:\Z^N\curvearrowright\mathcal{Z}$ be 
a strongly outer action of $\Z^N$. 
Let $\gamma'$ be the $\Z^{N-1}$-action 
generated by the first $N{-}1$ generators of $\gamma$. 
Let $\tilde\gamma_{\xi_N}
\in\Aut_{\T^{N-1}}(\mathcal{Z}\rtimes_{\gamma'}\Z^{N-1})$ be 
the natural extension of the automorphism $\gamma_{\xi_N}$ of $\mathcal{Z}$. 
As mentioned in the proof of Lemma \ref{KKcomputation}, 
the canonical inclusion of $C^*(\Z^{N-1})$ 
into $\mathcal{Z}\rtimes_{\gamma'}\Z^{N-1}$ gives a $KK$-equivalence. 
Therefore $K_i(\tilde\gamma_{\xi_N})=\id$ for $i=0,1$. 
By \cite[Theorem 5.8]{MS}, we get the conclusion. 
\end{proof}

\begin{lem}\label{admissible}
Let $\gamma:\Z^N\curvearrowright\mathcal{Z}$ be 
a strongly outer, approximately representable action of $\Z^N$ 
on the Jiang-Su algebra and 
let $A$ be a unital simple infinite dimensional AF algebra such that 
$K_0(A)$ has no infinitesimal. 
Then, for any $\id\otimes\gamma$-cocycle 
$(u_g)_{g\in\Z^N}$ in $A\otimes\mathcal{Z}$ and $\ep>0$, 
there exists a unitary $v\in A\otimes\mathcal{Z}$ such that 
\[
\lVert u_{\xi_i}-v(\id\otimes\gamma_{\xi_i})(v^*)\rVert<\ep
\]
for each $i=1,2,\dots,N$. 
\end{lem}
\begin{proof}
Set $B=(A\otimes\mathcal{Z})\rtimes_{\id\otimes\gamma}\Z^N$. 
By Lemma \ref{Z} and \cite[Definition 2.5]{MS}, 
$B\cong A\otimes(\mathcal{Z}\rtimes_\gamma\Z^N)$ belongs to $\mathcal{C}_0$. 
By \cite[Lemma 2.4]{MS}, $B$ has tracial rank zero. 
Since the $K$-groups of $B$ are torsion free, 
$B$ is a unital simple AT algebra with real rank zero. 

Let $\iota:C^*(\Z^N)\to B$ be the canonical embedding. 
Define a homomorphism $\iota_u:C^*(\Z^N)\to B$ 
by $\iota_u(\lambda_g)=u_g\lambda^{\id\otimes\gamma}_g$ for $g\in \Z^N$. 
By means of Lemma \ref{KKcomputation}, one has $KK(\iota)=KK(\iota_u)$. 
It is clear that for any $\tau\in T(B)$ and $g\in\Z^N$ 
\[
\tau(\iota(\lambda_g))=\tau(\lambda^{\id\otimes\gamma}_g)
=\begin{cases}0&g\neq0\\1&g=0\end{cases}
\]
and 
\[
\tau(\iota_u(\lambda_g))=\tau(u_g\lambda^{\id\otimes\gamma}_g)
=\begin{cases}0&g\neq0\\1&g=0.\end{cases}
\]
It follows that 
$\tau\circ\iota$ equals $\tau\circ\iota_u$ for any $\tau\in T(B)$. 
Hence, by \cite[Theorem 3.4]{L07Trans} or \cite[Theorem 4.8]{M10}, 
the two homomorphisms $\iota$ and $\iota_u$ are 
approximately unitarily equivalent. 
Since $\id\otimes\gamma$ is approximately representable, 
by \cite[Corollary 4.10]{IM}, we can conclude the proof. 
\end{proof}

The following is the main result of this paper. 

\begin{thm}\label{SCConUHF}
Let $A$ be a UHF algebra of infinite type. 
Then any two strongly outer actions of $\Z^N$ on $A$ are 
strongly cocycle conjugate to each other. 
\end{thm}
\begin{proof}
By Theorem \ref{CCofRohlin1} and Theorem \ref{Rohlintype1}, 
any two strongly outer actions of $\Z^N$ on $A$ are 
cocycle conjugate to each other. 
It follows from the lemma above that 
they are strongly cocycle conjugate to each other. 
\end{proof}

For certain simple AF algebras, 
we can strengthen Theorem \ref{CCofasymp} as follows. 

\begin{thm}\label{SCCofasymp}
Let $A$ be a unital simple AF algebra 
such that $K_0(A)$ has no infinitesimal. 
Suppose that $A$ has finitely many extremal tracial states. 
Let $\alpha,\beta:\Z^N\curvearrowright A$ be 
strongly outer, asymptotically representable actions of $\Z^N$. 
Then they are strongly cocycle conjugate. 
\end{thm}
\begin{proof}
This is an immediate consequence of 
Theorem \ref{CCofasymp} and Lemma \ref{admissible}. 
\end{proof}

\end{document}